\documentstyle[11pt,amssymb,amsfonts]{article}

\begin{document}
\newcommand{\p}{\parallel }
\makeatletter \makeatother
\newtheorem{th}{Theorem}[section]
\newtheorem{lem}{Lemma}[section]
\newtheorem{de}{Definition}[section]
\newtheorem{rem}{Remark}[section]
\newtheorem{cor}{Corollary}[section]
\renewcommand{\theequation}{\thesection.\arabic {equation}}

\title{{\bf Lower dimensional volumes and the Kastler-Kalau-Walze type
theorem for Manifolds with Boundary} }
\author{  Yong Wang \\}

\date{}
\maketitle

\begin{abstract} In this paper, we define lower dimensional
volumes of spin manifolds with boundary. We compute the lower
dimensional volume ${\rm Vol}^{(2,2)}$ for $5$-dimensional and
$6$-dimensional spin manifolds with boundary and we also get the
Kastler-Kalau-Walze type theorem in this case.\\

\noindent{\bf MSC:}\quad 58G20; 53A30; 46L87\\
 \noindent{\bf Keywords:}\quad
Lower dimensional volumes, noncommutative residue for manifolds with
boundary; gravitational action for manifolds with
boundary\\

\end{abstract}

\section{Introduction}
\quad The noncommutative residue found in [Gu] and [Wo] plays a
prominent role in noncommutative geometry. In [C1], Connes used the
noncommutative residue to derive a conformal 4-dimensional Polyakov
action analogy. In [C2], Connes proved that the noncommutative
residue on a compact manifold $M$ coincided with the Dixmier's trace
on pseudodifferential operators of order $-{\rm {dim}}M$. Several
years ago, Connes made a challenging observation that the
noncommutative residue of the square of the inverse of the Dirac
operator was proportional to the Einstein-Hilbert action, which we
call the Kastler-Kalau-Walze theorem. In [K], Kastler gave a
brute-force proof of this theorem. In [KW], Kalau and Walze proved
this theorem in the normal coordinates system simultaneously. In
[A], Ackermann gave a note on a new proof of this theorem
by means of the heat kernel expansion.\\
\indent On the other hand, Fedosov et al defined a noncommutative
residue on Boutet de Monvel's algebra and proved that it was a
unique continuous trace in [FGLS]. In [Wa1] and [Wa2], we
generalized some results in [C1] and [U] to the case of manifolds
with boundary . In [Wa3], We proved a Kastler-Kalau-Walze type
theorem for the Dirac operator and the signature operator for
$3,4$-dimensional manifolds with boundary. Recently, Ponge defined
lower dimensional volumes of Riemannian manifolds by the Wodzicki
residue in [Po].  The motivation of this paper is to find a
Kastler-Kalau-Walze type theorem for higher dimensional manifolds
with boundary and generalize the definition of lower dimensional volumes to manifolds with boundary.\\
\indent This paper is organized as follows: In Section 2, we define
lower dimensional volumes of spin manifolds with boundary. In
Section $3$, for $6$-dimensional spin manifolds with boundary and
the associated Dirac operator $D$, we compute the lower dimensional
volume ${\rm Vol}^{(2,2)}_6$ and get a Kastler-Kalau-Walze type
theorem in this case. In Section 4, when $\partial M$ is flat, we
can define $\int_{\partial M }{\rm res}_{2,2}(D^{-2},D^{-2})$ and
$\int_{\partial M} {\rm res}_{2,3}(D^{-2},D^{-2})$ (see Section 3)
and get that the gravitational action for $\partial M$ is
proportional to $\int_{\partial M} {\rm res}_{2,2}(D^{-2},D^{-2})$
and $\int_{\partial M }{\rm res}_{2,3}(D^{-2},D^{-2})$, which gives
two kinds of operator theoretic explanations of the gravitational
action for boundary. For $5$-dimensional spin manifolds with
boundary and the associated Dirac operator $D$, we compute the lower
dimensional
volume ${\rm Vol}^{(2,2)}_5$.\\

\section{Lower dimensional
volumes of spin manifolds with boundary}

   \quad In order to define lower dimensional volumes of spin manifolds with
boundary, we need some basic facts and formulae about Boutet de
Monvel's calculus and the definition of noncommutative residue for
manifolds with boundary. We can find them in Section 2,3 in [Wa1]
and Section 2.1 in [Wa3].\\
  \indent Let $M$ be a $n$-dimensional compact oriented spin manifold
with boundary $\partial M$. We assume that the metric $g^M$ on $M$
has the following form near the boundary,
$$ g^M=\frac{1}{h(x_n)}g^{\partial M}+dx_n^2,\eqno(2.1)$$
where $g^{\partial M}$ is the metric on  ${\partial M}$. $h(x_n)\in
C^{\infty}([0,1))=\{\widetilde{h}|_{[0,1)}|\widetilde{h}\in
C^{\infty}((-\varepsilon,1))\}$ for some $\varepsilon>0$ and
satisfies $h(x_n)>0,~h(0)=1$ where $x_n$ denotes the normal
directional coordinate. Let $U\subset M$ be a collar neighborhood of
$\partial M$ which is diffeomorphic to $\partial M\times [0,1)$. By
the definition of $C^{\infty}([0,1))$ and $h>0$, there exists
$\widetilde{h}\in C^{\infty}((-\varepsilon,1))$ such that
$\widetilde{h}|_{[0,1)}=h$ and $\widetilde{h}>0$ for some
sufficiently small $\varepsilon>0$. Then there exists a metric
$\widehat{g}$ on $\widehat{M}=M\cup_{\partial M}\partial M\times
(-\varepsilon,0]$ which has the form on $U\cup_{\partial M}\partial
M\times (-\varepsilon,0]$
$$ \widehat{g}=\frac{1}{\widetilde{h}(x_n)}g^{\partial M}+dx_n^2,\eqno(2.2)$$
such that $\widehat{g}|_M=g.$ We fix a metric $\widehat{g}$ on the
$\widehat{M}$ such that $\widehat{g}|_M=g$. We can get the spin
structure on $\widehat{M}$ by extending the spin structure on $M.$
Let $D$ be the Dirac operator associated to $\widehat{g}$ on the
spinors bundle $S(T\widehat{M})$. Let $p_1,p_2$ be nonnegative
integers
and $p_1+p_2\leq n$.\\

\noindent {\bf Definition 2.1} Lower dimensional volumes of spin
manifolds with boundary are defined by ${\rm Vol}^{(p_1,p_2)}_nM:=
\widetilde{{\rm Wres}}[\pi^+D^{-p_1}\circ\pi^+D^{-p_2}]$ (for the
related
definitions, see [Wa1], Section 2, 3).\\

  \indent Denote by $\sigma_l(A)$ the $l$-order
symbol of an operator $A$. By (2.1.4)-(2.1.8) in [Wa3], we get
$$\widetilde{{\rm Wres}}[\pi^+D^{-p_1}\circ\pi^+D^{-p_2}]=\int_M\int_{|\xi|=1}{\rm
trace}_{S(TM)}[\sigma_{-n}(D^{-p_1-p_2})]\sigma(\xi)dx+\int_{\partial
M}\Phi,\eqno(2.3)$$ \noindent where
$$\Phi=\int_{|\xi'|=1}\int^{+\infty}_{-\infty}\sum^{\infty}_{j, k=0}
\sum\frac{(-i)^{|\alpha|+j+k+1}}{\alpha!(j+k+1)!}$$
$$\times {\rm trace}_{S(TM)}
[\partial^j_{x_n}\partial^\alpha_{\xi'}\partial^k_{\xi_n}
\sigma^+_{r}(D^{-p_1})(x',0,\xi',\xi_n)\times
\partial^\alpha_{x'}\partial^{j+1}_{\xi_n}\partial^k_{x_n}\sigma_{l}
(D^{-p_2})(x',0,\xi',\xi_n)]d\xi_n\sigma(\xi')dx',\eqno(2.4)$$
\noindent where the sum is taken over $
r-k-|\alpha|+l-j-1=-n,~~r\leq -p_1,l\leq -p_2$.
 Since
$[\sigma_{-n}(D^{-p_1-p_2})]|_M$ has the same expression as
$\sigma_{-n}(D^{-p_1-p_2})$ in the case of manifolds without
boundary, so locally we can use the computations in [K], [KW], [Po]
to compute the first term. The following proposition is the
motivation of the definition of lower dimensional volumes of spin
manifolds with boundary.\\

\noindent {\bf Proposition 2.2} 1){\it When $p_1+p_2=n$, then ${\rm
Vol}^{(p_1,p_2)}_nM=c_0{\rm Vol}_M.$}\\
2) when $p_1+p_2\equiv n {\rm mod} 1$, ${\rm
Vol}^{(p_1,p_2)}_nM=\int_{\partial M}\Phi.$\\
3)$${\rm Vol}^{(1,1)}_4=-\frac{\Omega_4}{3}\int_Ms{\rm dvol}_M;~~
{\rm Vol}^{(1,1)}_3=c_1{\rm Vol}_{\partial M}\eqno(2.5)$$ {\it where
$c_0,c_1$ are constants and $s$ is the scalar curvature.}\\

\noindent {\bf Proof.} 1) comes from (2.4) and (2.2) in [Po]. 2)
comes from the proposition 2.3 and 3.2 in [Po]. 3) comes from
Theorems 2.5 and
5.1 in [Wa3]. $\Box$\\

\section{A Kastler-Kalau-Walze type theorem for $6$-dimensional spin
manifolds with boundary }

\quad In this section, We compute the lower dimensional volume ${\rm
Vol}^{(2,2)}_6$ for $6$-dimensional spin manifolds with
boundary and get a Kastler-Kalau-Walze type theorem in this case.\\
\indent Firstly, we recall the symbol expansion of $D^{-2}$ in [Ka].
 Recall the
definition of the Dirac operator $D$ (see [BGV], [Y]). Let
$\nabla^L$ denote the Levi-civita connection about $g^M$.
 In the local coordinates $\{x_i; 1\leq i\leq n\}$ and
 the fixed orthonormal frame $\{\widetilde{e_1},\cdots,\widetilde{e_n}\}$, the connection matrix $(\omega_{s,t})$
is defined by
$$\nabla^L(\widetilde{e_1},\cdots,\widetilde{e_n})= (\widetilde{e_1},\cdots,\widetilde{e_n})(\omega_{s,t}).\eqno(3.1)$$
$c(\widetilde{e_i})$ denotes the Clifford action. The Dirac operator
$$D=\sum^n_{i=1}c({\partial_i})(\partial_i+\delta_i);~
\delta_i=-\frac{1}{4}\sum_{s,t}\omega_{s,t}(\partial_i)c(\widetilde{e_s})c(\widetilde{e_t}).\eqno(3.2)$$
Let $g^{ij}=g(dx_i,dx_j)$ and
$$\nabla^L_{\partial_i}\partial_j=\sum_k\Gamma_{ij}^k\partial_k;
~\Gamma^k=g^{ij}\Gamma_{ij}^k, ~\delta^j=g^{ij}\delta_i.\eqno(3.3)$$
Let the cotangent vector $\xi=\sum \xi_jdx_j$ and
$\xi^j=g^{ij}\xi_i$. Then we have\\

\noindent {\bf Lemma 3.1}([Ka]) $$\sigma_{-2}(D^{-2})=|\xi|^{-2};~
\sigma_{-3}(D^{-2})=-\sqrt{-1}|\xi|^{-4}\xi_k(\Gamma^k-2\delta^k)-\sqrt{-1}|\xi|^{-6}2\xi^j\xi_\alpha\xi_\beta
\partial_jg^{\alpha\beta}.\eqno(3.4)$$

\indent Since $\Phi$ is a global form on $\partial M$, so for any
fixed point $x_0\in\partial M$, we can choose the normal coordinates
$U$ of $x_0$ in $\partial M$ (not in $M$) and compute $\Phi(x_0)$ in
the coordinates $\widetilde{U}=U\times [0,1)\subset M$ and the
metric $\frac{1}{h(x_n)}g^{\partial M}+dx_n^2.$ For details, see
Section 2.2.2 in [Wa3].\\
\indent Now we can compute $\Phi$ (see formula (2.4) for the
definition of $\Phi$), since the sum is taken over $
-r-l+k+j+|\alpha|=-5,~~r,l\leq-2,$ then we have the following five cases:\\

\noindent  {\bf case a)~I)}~$r=-2,~l=-2~k=j=0,~|\alpha|=1$\\

\noindent By (2.4), we get
$${\rm case~a)~I)}=-\int_{|\xi'|=1}\int^{+\infty}_{-\infty}\sum_{|\alpha|=1}
{\rm trace}
[\partial^\alpha_{\xi'}\pi^+_{\xi_n}\sigma_{-2}(D^{-2})\times
\partial^\alpha_{x'}\partial_{\xi_n}\sigma_{-2}(D^{-2})](x_0)d\xi_n\sigma(\xi')dx',\eqno(3.5)$$
By Lemma 2.2 in [Wa3], for $i<n$, then
$$\partial_{x_i}\sigma_{-2}(D^{-2})(x_0)=\partial_{x_i}{(|\xi|^{-2})}(x_0)=
-\frac{\partial_{x_i}(|\xi|^2)(x_0)}{|\xi|^4}=0,\eqno(3.6)$$
\noindent so case a) I) vanishes.\\

\noindent  {\bf case a)~II)}~$r=-2,~l=-2~k=|\alpha|=0,~j=1$\\

\noindent By (2.4), we get
$${\rm case~
a)~II)}=-\frac{1}{2}\int_{|\xi'|=1}\int^{+\infty}_{-\infty} {\rm
trace} [\partial_{x_n}\pi^+_{\xi_n}\sigma_{-2}(D^{-2})\times
\partial_{\xi_n}^2\sigma_{-2}(D^{-2})](x_0)d\xi_n\sigma(\xi')dx',\eqno(3.7)$$
\noindent By Lemma 2.2 in [Wa3], we have\\
$$\partial_{x_n}\sigma_{-2}(D^{-2})(x_0)|_{|\xi'|=1}=-\frac{h'(0)}{(1+\xi_n^2)^2}\eqno(3.8)$$
By (3.8) and the Cauchy integral formula and (2.1.1) in [Wa3], then\\
\begin{eqnarray*}
\pi^+_{\xi_n}\partial_{x_n}\sigma_{-2}(D^{-2})(x_0)|_{|\xi'|=1}
&=&-h'(0)\frac{1}{2\pi i}{\rm lim}_{u\rightarrow
0^-}\int_{\Gamma^+}\frac{\frac{1}{(\eta_n+i)^2(\xi_n+iu-\eta_n)}}
{(\eta_n-i)^2}d\eta_n\\
&=&h'(0)\frac{i\xi_n+2}{4(\xi_n-i)^2},~~~~~~~~~~~~~~~~~~~~~~~~~~~~~~~~~~~~(3.9)
\end{eqnarray*}
$$\partial^2_{\xi_n}(|\xi|^{-2})(x_0)=\frac{-2+6\xi_n^2}{(1+\xi_n^2)^3}.\eqno(3.10)$$
We note that
\begin{eqnarray*}
\int_{-\infty}^{\infty}\frac{i\xi_n+2}{(\xi_n-i)^2}\times
\frac{-1+3\xi_n^2}{(1+\xi_n^2)^3}d\xi_n
&=&\int_{\Gamma^+}\frac{3i\xi_n^3+6\xi_n^2-i\xi_n-2}{(\xi_n-i)^5(\xi_n+i)^3}d\xi_n\\
&=&\frac{2\pi
i}{4!}\left[\frac{3i\xi_n^3+6\xi_n^2-i\xi_n-2}{(\xi_n+i)^3}\right]^{(4)}|_{\xi_n=i}=\frac{5\pi}{16},~~~(3.11)
\end{eqnarray*}
Since $n=6$, ${\rm tr}_{S(TM)}[{\rm id}]={\rm dim}(\wedge^*(3))=8.$
So by (3.7)-(3.11), we get {\bf case a} II)$=-\frac{5}{8}\pi
h'(0)\Omega_4dx',$ where $\Omega_4$ is the
canonical volume of $S^4$.\\

\noindent  {\bf case a)~III)}~$r=-2,~l=-2~j=|\alpha|=0,~k=1$\\

\noindent By (2.4) and an integration by parts, we get
$${\rm case~ a)~III)}=-\frac{1}{2}\int_{|\xi'|=1}\int^{+\infty}_{-\infty}
{\rm trace} [\partial_{\xi_n}\pi^+_{\xi_n}\sigma_{-2}(D^{-2})\times
\partial_{\xi_n}\partial_{x_n}\sigma_{-2}(D^{-2})](x_0)d\xi_n\sigma(\xi')dx'$$
$$=\frac{1}{2}\int_{|\xi'|=1}\int^{+\infty}_{-\infty} {\rm trace}
[\partial_{\xi_n}^2\pi^+_{\xi_n}\sigma_{-2}(D^{-2})\times
\partial_{x_n}\sigma_{-2}(D^{-2})](x_0)d\xi_n\sigma(\xi')dx'.\eqno(3.12)$$
\noindent By Lemma 2.2 in [Wa3] , we have\\
$$\partial_{\xi_n}^2\pi_{\xi_n}^+\sigma_{-2}(D^{-2})(x_0)|_{|\xi'|=1}=\frac{-i}{(\xi_n-i)^3}.\eqno(3.13)$$
By (3.8) and (3.13), we have $$ {\rm {\bf case~
a)~III)}}=4ih'(0)\int_{|\xi'|=1}\int^{+\infty}_{-\infty}\int_{\Gamma^+}\frac{1}{(\xi_n-i)^5(\xi_n+i)^2}d\xi_n\sigma(\xi')dx'
=
\frac{5}{8}\pi h'(0)\Omega_4dx'.$$ Thus the sum of {\bf case~
a)~II)} and {\bf case~ a)~III)} is zero.\\

\noindent  {\bf case b)}~$r=-2,~l=-3,~k=j=|\alpha|=0$\\

\noindent By (2.4) and an integration by parts, we get
$${\rm case~ b)}=-i\int_{|\xi'|=1}\int^{+\infty}_{-\infty}
{\rm trace} [\pi^+_{\xi_n}\sigma_{-2}(D^{-2})\times
\partial_{\xi_n}\sigma_{-3}(D^{-2})](x_0)d\xi_n\sigma(\xi')dx'$$
$$=i\int_{|\xi'|=1}\int^{+\infty}_{-\infty}
{\rm trace} [\partial_{\xi_n}\pi^+_{\xi_n}\sigma_{-2}(D^{-2})\times
\sigma_{-3}(D^{-2})](x_0)d\xi_n\sigma(\xi')dx'.\eqno(3.14)$$
\noindent By Lemma 2.2 in [Wa3], we have\\
$$\partial_{\xi_n}\pi_{\xi_n}^+\sigma_{-2}(D^{-2})(x_0)|_{|\xi'|=1}=\frac{i}{2(\xi_n-i)^2}.\eqno(3.15)$$
In the normal coordinate, $g^{ij}(x_0)=\delta_i^j$ and
$\partial_{x_j}(g^{\alpha\beta})(x_0)=0,$ {\rm if
}$j<n;~=h'(0)\delta^\alpha_\beta,~{\rm if }~j=n.$ So by Lemma A.2 in
[Wa3], we have $\Gamma^n(x_0)=\frac{5}{2}h'(0)$ and
$\Gamma^k(x_0)=0$ for $k<n$. By the definition of $\delta^k$ and
Lemma 2.3 in [Wa3], we have $\delta^n(x_0)=0$ and
$\delta^k=\frac{1}{4}h'(0)c(\widetilde{e_k})c(\widetilde{e_n})$ for
$k<n$. So
$$\sigma_{-3}(D^{-2})(x_0)|_{|\xi'|=1}
=-\sqrt{-1}|\xi|^{-4}\xi_k(\Gamma^k-2\delta^k)(x_0)|_{|\xi'|=1}-\sqrt{-1}|\xi|^{-6}2\xi^j\xi_\alpha\xi_\beta
\partial_jg^{\alpha\beta}(x_0)|_{|\xi'|=1}$$
$$=\frac{-i}{(1+\xi_n^2)^2}(-\frac{1}{2}h'(0)\sum_{k<n}\xi_k
c(\widetilde{e_k})c(\widetilde{e_n})+\frac{5}{2}h'(0)\xi_n)-\frac{2ih'(0)\xi_n}{(1+\xi_n^2)^3}.\eqno(3.16)$$
We note that $\int_{|\xi'|=1}\xi_1\cdots\xi_{2q+1}\sigma(\xi')=0$,
so the first term in (3.16) has no contribution for computing case
b).\\
\begin{eqnarray*}
{\bf case~ b)}&=& i h'(0)\int_{|\xi'|=1}\int^{+\infty}_{-\infty}
{\rm trace}
\left[\frac{1}{2(\xi_n-i)^2}\times\left(\frac{\frac{5}{2}\xi_n}{(1+\xi_n^2)^2}
+\frac{2\xi_n}{(1+\xi_n^2)^3}\right)\right]d\xi_n\sigma(\xi')dx'\\
&=&2ih'(0)\Omega_4\int_{\Gamma^+}\frac{5\xi_n^3+9\xi_n}{(\xi_n-i)^5(\xi_n+i)^3}d\xi_ndx'\\
&=&2ih'(0)\Omega_4\frac{2\pi
i}{4!}\left[\frac{5\xi_n^3+9\xi_n}{(\xi_n+i)^3}\right]^{(4)}|_{\xi_n=i}dx'
=-\frac{15}{8}\pi h'(0)\Omega_4dx'.~~~~~~~~~(3.17)
\end{eqnarray*}

\noindent {\bf  case c)}~$r=-3,~l=-2,~k=j=|\alpha|=0$\\

\noindent By (2.4), we get
$${\rm case~ c)}=-i\int_{|\xi'|=1}\int^{+\infty}_{-\infty}
{\rm trace} [\pi^+_{\xi_n}\sigma_{-3}(D^{-2})\times
\partial_{\xi_n}\sigma_{-2}(D^{-2})](x_0)d\xi_n\sigma(\xi')dx'.\eqno(3.18)$$
By the Leibniz rule, trace property and "++" and "-~-" vanishing
after the integration over $\xi_n$ (for details, see [FGLS]), then
\begin{eqnarray*}
& &\int^{+\infty}_{-\infty}{\rm trace}
[\pi^+_{\xi_n}\sigma_{-3}(D^{-2})\times
\partial_{\xi_n}\sigma_{-2}(D^{-2})]d\xi_n\\
&=&\int^{+\infty}_{-\infty}{\rm tr} [\sigma_{-3}(D^{-2})\times
\partial_{\xi_n}\sigma_{-2}(D^{-2})]d\xi_n
-\int^{+\infty}_{-\infty}{\rm tr}
[\pi^-_{\xi_n}\sigma_{-3}(D^{-2})\times
\partial_{\xi_n}\sigma_{-2}(D^{-2})]d\xi_n\\
&=&\int^{+\infty}_{-\infty}{\rm tr} [\sigma_{-3}(D^{-2})\times
\partial_{\xi_n}\sigma_{-2}(D^{-2})]d\xi_n-\int^{+\infty}_{-\infty}{\rm
tr} [\pi^-_{\xi_n}\sigma_{-3}(D^{-2})\times
\partial_{\xi_n}\pi^+_{\xi_n}\sigma_{-2}(D^{-2})]d\xi_n\\
&=&\int^{+\infty}_{-\infty}{\rm tr} [\sigma_{-3}(D^{-2})\times
\partial_{\xi_n}\sigma_{-2}(D^{-2})]d\xi_n-\int^{+\infty}_{-\infty}{\rm
tr} [\sigma_{-3}(D^{-2})\times
\partial_{\xi_n}\pi^+_{\xi_n}\sigma_{-2}(D^{-2})]d\xi_n\\
&=&\int^{+\infty}_{-\infty}{\rm tr} [
\partial_{\xi_n}\sigma_{-2}(D^{-2})\times\sigma_{-3}(D^{-2})]d\xi_n+\int^{+\infty}_{-\infty}{\rm
tr}[\partial_{\xi_n}\sigma_{-3}(D^{-2})\times
\pi^+_{\xi_n}\sigma_{-2}(D^{-2})]d\xi_n\\
&=&\int^{+\infty}_{-\infty}{\rm tr} [
\partial_{\xi_n}\sigma_{-2}(D^{-2})\times\sigma_{-3}(D^{-2})]d\xi_n+\int^{+\infty}_{-\infty}{\rm
tr}[\pi^+_{\xi_n}\sigma_{-2}(D^{-2})\times
\partial_{\xi_n}\sigma_{-3}(D^{-2})]d\xi_n ~~~(3.19)
\end{eqnarray*}
By (3.19), we have
$${\rm {\bf case~ c)}}={\rm {\bf case~ b)}}-i\int_{|\xi'|=1}\int^{+\infty}_{-\infty}{\rm
tr}[\partial_{\xi_n}\sigma_{-2}(D^{-2})\times
\sigma_{-3}(D^{-2})]d\xi_n\sigma(\xi')dx'.\eqno(3.20)$$ We note that
we can not get the sum of {\bf case b)} and {\bf case c)} is zero by
computations in (3.19). In order to compute {\bf case c)}, we only
need compute the last term in (3.20). By (3.16) and
$$\partial_{\xi_n}\sigma_{-2}(D^{-2})(x_0)|_{|\xi'|=1}=-\frac{2\xi_n}{(\xi_n^2+1)^2},\eqno(3.21)$$
we have $$-i\int_{|\xi'|=1}\int^{+\infty}_{-\infty}{\rm
tr}[\partial_{\xi_n}\sigma_{-2}(D^{-2})\times
\sigma_{-3}(D^{-2})]d\xi_n\sigma(\xi')dx'$$
$$=8h'(0)\int_{|\xi'|=1}\int^{+\infty}_{-\infty}\frac{5\xi_n^4+9\xi_n^2}{(\xi_n+i)^5(\xi_n-i)^5}d\xi_n\sigma(\xi')dx'
=\frac{15}{4}\pi h'(0)\Omega_4dx'.\eqno(3.22)$$
 By (3.17), (3.20) and (3.22), we have the sum of {\bf case b)} and {\bf case c)} is
 zero. Now $\Phi$ is the sum of the cases a), b) and c), so is zero. Then we get\\

\noindent {\bf Theorem 3.2}~~{\it  Let  $M$ be a $6$-dimensional
compact spin manifold with the boundary $\partial M$ and the metric
$g^M$ as above and $D$ be the Dirac operator on $\widehat{M}$, then}
$${\rm Vol}^{(2,2)}_6= \widetilde{{\rm
Wres}}[(\pi^+D^{-2})^2]=-\frac{5\Omega_6}{3}\int_Ms{\rm
dvol}_M.\eqno(3.23)$$

\section{The gravitational action for $6$-dimensional manifolds with boundary}
 \quad Firstly, we recall the Einstein-Hilbert action for manifolds with boundary (see [H] or [B]),
$$I_{\rm Gr}=\frac{1}{16\pi}\int_Ms{\rm dvol}_M+2\int_{\partial M}K{\rm dvol}_{\partial_M}:=I_{\rm {Gr,i}}+I_{\rm {Gr,b}},\eqno(4.1)$$
\noindent  where
$$ K=\sum_{1\leq i,j\leq {n-1}}K_{i,j}g_{\partial M}^{i,j};~~K_{i,j}=-\Gamma^n_{i,j},\eqno(4.2)$$
\noindent and $K_{i,j}$ is the second fundamental form, or extrinsic
curvature. Taking the metric in Section 2, then by Lemma A.2 in
[Wa3],
 for $n=6$, then
$$K(x_0)==-\frac{5}{2}h'(0);~
 I_{\rm {Gr,b}}=-5h'(0){\rm Vol}_{\partial M}.\eqno(4.3)$$
\indent Let $M$ be $6$-dimensional manifolds with boundary and
$P,P'$ be two pseudodifferential operators with transmission
property (see [Wa1] or [RS]) on $\widehat M$. Motivated by (2.4), we
define locally
 $${\rm
res}_{2,2}(P,P'):=-\frac{1}{2}\int_{|\xi'|=1}\int^{+\infty}_{-\infty}
{\rm trace} [\partial_{x_n}\pi^+_{\xi_n}\sigma_{-2}(P)\times
\partial_{\xi_n}^2\sigma_{-2}(P')]d\xi_n\sigma(\xi')dx';\eqno(4.4)$$
$${\rm res}_{2,3}(P,P'):=-i\int_{|\xi'|=1}\int^{+\infty}_{-\infty}
{\rm trace} [\pi^+_{\xi_n}\sigma_{-2}(P)\times
\partial_{\xi_n}\sigma_{-3}(P')]d\xi_n\sigma(\xi')dx'.\eqno(4.5)$$
\noindent By (4.4),(4.5), so
$$ {\rm case~ a)~ II)}={\rm res}_{2,2}(D^{-2},D^{-2});~{\rm case~ b)}={\rm res}_{2,3}(D^{-2},D^{-2}).\eqno(4.6)$$
\noindent Now, we assume $\partial M$ is flat, then
$\{dx_i=e_i\},~g^{\partial M}_{i,j}=\delta_{i,j},~\partial
_{x_s}g^{\partial M}_{i,j}=0$. So ${\rm res}_{2,2}(D^{-2},D^{-2})$
and ${\rm res}_{2,3}(D^{-2},D^{-2})$ are two global forms locally
defined by the aboved oriented orthonormal basis $\{dx_i\}$. By case
a) II) and case b),
then we have:\\

\noindent {\bf Theorem 4.1}~~{\it  Let  $M$ be a $6$-dimensional
compact spin manifold with the boundary $\partial M$ and the metric
$g^M$ as above and $D$ be the Dirac operator on $\widehat{M}$.
Assume $\partial M$ is flat, then}
$$\int_{\partial M}{\rm res}_{2,2}(D^{-2},D^{-2})=\frac{\pi}{8}\Omega_4I_{\rm {Gr,b}};\eqno(4.7)$$
$$\int_{\partial M}{\rm res}_{2,3}(D^{-2},D^{-2})=\frac{3\pi}{8}\Omega_4I_{\rm {Gr,b}}.\eqno(4.8)$$\\

Nextly, for $5$-dimensional spin manifolds with boundary, we compute
${\rm Vol}^{(2,2)}_5$. By Proposition 2.2 (2), we have
$$\widetilde{{\rm
Wres}}[(\pi^+D^{-2})^2]=\int_{\partial M}\Phi.\eqno(4.9)$$ \noindent
When $n=5$, then in (2.4), $ r-k-|\alpha|+l-j-1=-5,~~r,l\leq-2$, so
we get $r=l=-2,~k=|\alpha|=j=0,$ then
$$\Phi=\int_{|\xi'|=1}\int^{+\infty}_{-\infty}
 {\rm trace}_{S(TM)}
[ \sigma^+_{-2}(D^{-2})(x',0,\xi',\xi_n)\times
\partial_{\xi_n}\sigma_{-2}
(D^{-2})(x',0,\xi',\xi_n)]d\xi_5\sigma(\xi')dx'.\eqno(4.10)$$ By
(3.21) and
$\pi^+_{\xi_n}\sigma_{-2}(x_0)|_{|\xi'|=1}=\frac{1}{2i(\xi_n-i)}$
and ${\rm tr}(id)={\rm dim}(S(TM))=4$, we can get ${\rm
Vol}^{(2,2)}_5=\frac{\pi i}{2}\Omega_3{\rm Vol}_{\partial M}.$  By
$I_{\rm
{Gr,b}}=-4h'(0){\rm Vol}_{\partial M}$, we have\\

\noindent {\bf Theorem 4.2}~~{\it Let $M$ be a $5$-dimensional
compact spin manifold with the boundary $\partial M$ and the metric
$g^M$ as in Section 2 and $D$ be the Dirac operator on $\widehat{M}$
, then}
$$ {\rm
Vol}^{(2,2)}_5=\widetilde{{\rm Wres}}[(\pi^+D^{-2})^2]=\frac{\pi
i}{2}\Omega_3{\rm Vol}_{\partial M},\eqno(4.11)$$ $$I_{\rm
{Gr,b}}=\frac{8ih'(0)}{\pi\Omega_3}\widetilde{{\rm
Wres}}[(\pi^+D^{-2})^2],\eqno(4.12)$${\it where ${\rm Vol}_{\partial
M}$ denotes the canonical volume of
${\partial M}.$}\\

 \noindent {\bf Remark 4.3}~By Theorem 4.2, we know
that $\widetilde{{\rm Wres}}[(\pi^+D^{-1})^2]$ is proportional to
the gravitational action for boundary for $5$-dimensional manifolds
with boundary. But the constant depends on $h'(0)$.\\

\noindent{\bf Acknowledgement:}~~The author is indebted to Professor
Weiping Zhang for the suggestion to consider this problem. This work
was supported by NSFC No.10801027.\\ \\

\noindent{\bf References}\\

\noindent [A] T. Ackermann, {\it A note on the Wodzicki residue,} J.
Geom. Phys., 20, 404-406, 1996.\\
\noindent [B]  N. H. Barth, {\it The fourth-order gravitational
action for manifolds with boundaries,}
Class. Quantum Grav. 2, 497-513, 1985.\\
\noindent [BGV] N. Berline, E. Getzler, M. Vergne, {\it Heat Kernals
and Dirac Operators,}
Springer-Verlag, Berlin, 1992.\\
\noindent [C1] A. Connes, {\it Quantized calculus and applications,}
XIth International Congress of Mathematical Physics (paris,1994),
15-36, Internat Press, Cambridge, MA, 1995.\\
\noindent [C2] A. Connes. {\it The action functinal in
noncommutative
geometry,} Comm. Math. Phys., 117:673-683, 1998.\\
\noindent [FGLS] B. V. Fedosov, F. Golse, E. Leichtnam, and E.
Schrohe. {\it The noncommutative residue for manifolds with
boundary,} J. Funct.
Anal, 142:1-31,1996.\\
 \noindent [Gu] V.W. Guillemin, {\it A new proof of Weyl's
formula on the asymptotic distribution of eigenvalues}, Adv. Math.
55 no.2, 131-160, 1985.\\
\noindent [H] S. W. Hawking, {\it General Relativity. An Einstein
Centenary Survey,} Edited by S. W. Hawking and W.
Israel, Cambridge University Press,Cambridge-New York, 1979.\\
\noindent [K] D. Kastler, {\it The Dirac operator and gravitiation,}
Commun. Math. Phys, 166:633-643, 1995.\\
\noindent [KW] W. Kalau and M.Walze, {\it Gravity, non-commutative
geometry, and the Wodzicki residue,} J. Geom. Phys., 16:327-344, 1995.\\
\noindent[Po] R. Ponge, {\it Noncommutative geometry and lower
dimensional volumes in Riemannian geometry,} Lett. Math. Phys. 83
(2008), no. 1, 19--32.\\
 \noindent [RS] S. Rempel and B. W. Schulze,
Index theory of elliptic
boundary problems, Akademieverlag, Berlin, 1982.\\
\noindent [U] W. J. Ugalde, {\it Differential forms and the Wodzicki
residue,} J. Geom. Phys. 58 (2008), no. 12, 1739--1751.\\
\noindent [Wa1] Y. Wang, {\it Differential forms and the Wodzicki
residue for manifolds with boundary,} J. Geom. Phys.,
56:731-753, 2006.\\
 \noindent
[Wa2] Y. Wang, {\it Differential forms and the noncommutative
residue for manifolds with boundary in the non-product Case,} Lett.
math. Phys., 77:41-51, 2006.\\
\noindent [Wa3] Y. Wang, {\it Gravity and the noncommutative residue
for manifolds with boundary,} Lett. Math. Phys. 80 (2007), no. 1,
37--56.\\
 \noindent [Wo] M. Wodzicki,  {\it Local invariants of
spectral
asymmetry}, Invent.Math. 75 no.1 143-178, 1984.\\
\noindent [Y] Y. Yu, {\it The Index Theorem and The Heat Equation Method},
Nankai Tracts in Mathematics - Vol. 2, World Scientific Publishing, 2001.\\

 \indent{  School of Mathematics and Statistics,
Northeast Normal University, Changchun Jilin, 130024 , China }\\
\indent E-mail: {\it wangy581@nenu.edu.cn; }\\

\end{document}